\documentclass{amsart}
\usepackage{graphicx}
\usepackage{amsmath}
\usepackage{amssymb}
\usepackage{amsthm}
\usepackage{esint}

\newtheorem{theorem}{Theorem}[section]

\theoremstyle{Corollary}

\numberwithin{equation}{section}

\begin{document}

\title{A note on compact CR Yamabe soliton}

\author{Pak Tung Ho}
\address{Department of Mathematics, Sogang University, Seoul
121-742, Korea}

\email{ptho@sogang.ac.kr, paktungho@yahoo.com.hk}

\subjclass[2000]{Primary 53C44, 32V20; Secondary  53C21}

\date{August 13, 2014.}

\keywords{CR Yamabe flow; CR Yamabe soliton; CR manifold}

\begin{abstract}
In this paper, we show that the Webster scalar curvature of
any compact CR Yamabe soliton must be constant.

\end{abstract}

\maketitle

\section{Introduction}

Given a compact Riemannian manifold $(M,g)$, the Yamabe problem is to find
a metric conformal to $g$ such that it has constant scalar curvature. This was
solved by Aubin, Schoen, and Trudinger in \cite{Aubin0,Schoen,Trudinger}.
The (unnormalized) Yamabe flow was introduced to study the Yamabe problem, which is defined as
follows:
\begin{equation}\label{1}
\frac{\partial}{\partial t}g(t)=-R_{g(t)}g(t)\mbox{ for }t\geq 0,\hspace{2mm}g(0)=g.
\end{equation}
Here $R_{g(t)}$ is the scalar curvature of $g(t)$.
The existence and convergence of the Yamabe flow
have been studied in \cite{Brendle4,Brendle5,Chow,Schwetlick&Struwe,Ye}.
Yamabe soliton is a self-similar solution to the Yamabe flow. More precisely,
$g(t)$ is called a Yamabe soliton if there exist  a smooth function $\sigma(t)$ and
a $1$-parameter family of diffeomorphisms $\{\psi_t\}$ of $M$ such that
\begin{equation}\label{3}
g(t)=\sigma(t)\psi^*_t(g)
\end{equation}
is the solution of the Yamabe flow (\ref{1}), with $\sigma(0)=1$ and $\psi_0=id_M$.
The following is an alternative definition:
$(M,g)$ is called Yamabe soliton if there exist a
vector field $X$ and a constant $\rho\in\mathbb{R}$ such that
\begin{equation*}
(R_g-\rho)g=\mathcal{L}_Xg.
\end{equation*}
Here $R_g$ is the  scalar curvature of the metric $g$, and $\mathcal{L}_X$ is the Lie derivative in the direction of $X$.
Note that  these two  definitions are equivalent (see \cite{DiCerbo} for the proof).
Yamabe soliton has been studied by many authors. See \cite{Calvaruso,Cao,Chu,Daskalopoulos,DiCerbo,Hsu,Hsu1,Ma&Cheng,Ma&Miquel}
and the references therein. In particular, we mention the following theorem related to the main result in this paper,
which was obtained independently by di Cerbo and Disconzi in
\cite{DiCerbo} and by Hsu in \cite{Hsu}:
\begin{theorem}
Any compact Yamabe soliton must
have constant scalar curvature.
\end{theorem}

Suppose $(M,\theta)$ is a strictly pseudoconvex CR manifold of real dimension $2n+1$. The CR Yamabe problem is to
find a contact form conformal to $\theta$ such that it has constant Webster scalar curvature.
This was
solved by Jerison-Lee and Gamara-Yacoub in \cite{Gamara2,Gamara1,Jerison&Lee1,Jerison&Lee2,Jerison&Lee3}.
The (unnormalized) CR Yamabe flow is defined
as the evolution equation of the contact form $\theta(t)$:
\begin{equation}\label{2}
\frac{\partial}{\partial t}\theta(t)=-R_{\theta(t)}\,\theta(t)\mbox{ for }t\geq 0,\hspace{2mm}\theta(0)=\theta.
\end{equation}
Here $R_{\theta(t)}$ is the Webster scalar curvature of the contact form $\theta(t)$.
The CR Yamabe flow was introduced to tackle the CR Yamabe problem. See \cite{Chang&Chiu&Wu,Chang&Cheng,Ho2}
and the references therein. As in the Riemannian case,
CR Yamabe soliton is a self-similar solution to the CR Yamabe flow: we call $\theta(t)$ a CR Yamabe soliton
if there exist a smooth function $\sigma(t)$ and
a $1$-parameter family of CR diffeomorphisms $\{\psi_t\}$ of $M$ such that
\begin{equation}\label{5}
\theta(t)=\sigma(t)\psi^*_t(\theta)
\end{equation}
is the solution of the CR Yamabe flow (\ref{2}), with $\sigma(0)=1$ and $\psi_0=id_M$.

The following is our main result, which is the CR version of
the result of di Cerbo and Disconzi
in \cite{DiCerbo} and
 Hsu in \cite{Hsu} that we mentioned above.

\begin{theorem}\label{main}
If $(M,\theta(t))$ is a compact strictly pseudoconvex CR manifold satisfying (\ref{5}),
then the Webster scalar curvature of $(M,\theta(t))$ is constant.
\end{theorem}

\section{Proof}

In this section, we are going to prove Theorem \ref{main}.
We will consider the evolution of the quantity
\begin{equation}\label{2.0}
\frac{\int_MR_{\theta(t)}dV_{\theta(t)}}{(\int_MdV_{\theta(t)})^{\frac{n}{n+1}}}
\end{equation}
along the CR Yamabe flow (\ref{2}). Note that if
 $\theta(t)=u(t)^{\frac{2}{n}}\theta$ is the solution of the CR Yamabe flow (\ref{2}), then
$u(t)$ satisfies the following evolution equation:
\begin{equation}\label{2.1}
\frac{\partial}{\partial t}u(t)=-\frac{n}{2}R_{\theta(t)}u(t)\mbox{ for }t\geq 0.
\end{equation}
Therefore, by (\ref{2.1}), the volume form $dV_{\theta(t)}$ of $\theta(t)$ satisfies
\begin{equation}\label{2.2}
\frac{\partial}{\partial t}(dV_{\theta(t)})=\frac{\partial}{\partial t}(u(t)^{\frac{2n+2}{n}}dV_{\theta})=
\frac{2n+2}{n}u(t)^{\frac{2n+2}{n}-1}\frac{\partial u(t)}{\partial t}dV_{\theta}=-(n+1)R_{\theta(t)}dV_{\theta(t)},
\end{equation}
which implies that
\begin{equation}\label{2.3}
\frac{d}{dt}\left(\int_MdV_{\theta(t)}\right)=-(n+1)\int_MR_{\theta(t)}dV_{\theta(t)}.
\end{equation}
Since $\theta(t)=u(t)^{\frac{2}{n}}\theta$, $u(t)$ satisfies the CR Yamabe equation:
$$-(2+\frac{2}{n})\Delta_{\theta}u(t)+R_{\theta}u(t)=R_{\theta(t)}u(t)^{1+\frac{2}{n}}$$
where $\Delta_{\theta}$ is the sub-Laplacian of the contact form $\theta$. Differentiate it with respect to $t$,
one can derive that the following evolution equation of the
Webster scalar curvature $R_{\theta(t)}$ of $\theta(t)$:
(see \cite{Ho1} or \cite{Ho2} for the case of normalized CR Yamabe flow)
\begin{equation}\label{2.4}
\frac{\partial}{\partial t}R_{\theta(t)}=(n+1)\Delta_{\theta(t)} R_{\theta(t)}+R_{\theta(t)}^2.
\end{equation}
Here $\Delta_{\theta(t)}$ is the sub-Laplacian of the contact form $\theta(t)$.
Therefore, we have
\begin{equation}\label{2.5}
\begin{split}
&\frac{d}{dt}\left(\int_MR_{\theta(t)}dV_{\theta(t)}\right)\\
&=\int_M(\frac{\partial}{\partial t}R_{\theta(t)})dV_{\theta(t)}+\int_MR_{\theta(t)}\frac{\partial}{\partial t}(dV_{\theta(t)})\\
&=\int_M\Big((n+1)\Delta_{\theta(t)} R_{\theta(t)}+R_{\theta(t)}^2\Big)dV_{\theta(t)}-(n+1)\int_MR_{\theta(t)}^2dV_{\theta(t)}\\
&=-n\int_MR_{\theta(t)}^2dV_{\theta(t)}
\end{split}
\end{equation}
where we have used (\ref{2.2}) and (\ref{2.4}). Combining (\ref{2.3}) and (\ref{2.5}),
we obtain
\begin{equation}\label{2.6}
\begin{split}
\frac{d}{dt}\left(\frac{\int_MR_{\theta(t)}dV_{\theta(t)}}{(\int_MdV_{\theta(t)})^{\frac{n}{n+1}}}\right)
&=\frac{-n\left(\int_MR_{\theta(t)}^2dV_{\theta(t)}\right)\left(\int_MdV_{\theta(t)}\right)+n\left(\int_MR_{\theta(t)}dV_{\theta(t)}\right)^2}
{\left(\int_MdV_{\theta(t)}\right)^{\frac{n}{n+1}+1}}\leq 0
\end{split}
\end{equation}
where the last inequality follows from Cauchy-Schwarz inequality. This shows that
the quantity in (\ref{2.0}) is decreasing along the unnormalized CR Yamabe flow (\ref{2}).

On the other hand, the quantity in (\ref{2.0}) is invariant under the CR Yamabe soliton (\ref{5}).
To see this, note that if
$\theta(t)=\sigma(t)\psi^*_t(\theta)$ for some smooth function $\sigma(t)$ and
a $1$-parameter family of CR diffeomorphisms $\{\psi_t\}$ of $M$, then
$R_{\sigma(t)\psi^*_t(\theta)}=\sigma(t)^{-1}R_{\psi^*_t(\theta)}$
and $dV_{\sigma(t)\psi^*_t(\theta)}=\sigma(t)^{n+1}dV_{\psi^*_t(\theta)}$, which implies that
\begin{equation*}
\begin{split}
\frac{\int_MR_{\theta(t)}dV_{\theta(t)}}{(\int_MdV_{\theta(t)})^{\frac{n}{n+1}}}=
\frac{\sigma(t)^n\int_MR_{\psi^*_t(\theta)}dV_{\psi^*_t(\theta)}}{(\sigma(t)^{n+1}\int_MdV_{\psi^*_t(\theta)})^{\frac{n}{n+1}}}
=
\frac{\int_MR_{\theta}dV_{\theta}}{(\int_MdV_{\theta})^{\frac{n}{n+1}}}.
\end{split}
\end{equation*}
Therefore, we have
$$\frac{d}{dt}\left(\frac{\int_MR_{\theta(t)}dV_{\theta(t)}}{(\int_MdV_{\theta(t)})^{\frac{n}{n+1}}}\right)=0$$
under the CR Yamabe soliton (\ref{5}). This implies that
the inequality in (\ref{2.6}) is equality.
In particular,
$R_{\theta(t)}$ must be constant
by the equality case of the Cauchy-Schwarz inequality in (\ref{2.6}).
This completes the proof of
Theorem \ref{main}.

\bibliographystyle{amsplain}

\begin{thebibliography}{30}



\bibitem{Aubin0} T. Aubin, \'{E}quations diff\'{e}rentielles non lin\'{e}aires et probl\`{e}me de Yamabe concernant
la courbure scalaire. \textit{J. Math. Pures Appl. (9)} \textbf{55} (1976), 269--296.

\bibitem{Brendle4} S. Brendle, Convergence of the Yamabe flow for arbitrary initial
energy. \textit{J. Differential Geom.} \textbf{69} (2005), 217--278.

\bibitem{Brendle5} S. Brendle, Convergence of the Yamabe flow in dimension 6 and higher. \textit{Invent.
Math.} \textbf{170} (2007),  541--576.

\bibitem{Calvaruso}
G. Calvaruso
and A. Zaeim, A complete classification of Ricci and Yamabe solitons of non-reductive homogeneous 4-spaces.
\textit{J. Geom. Phys.} \textbf{80} (2014), 15--25.


\bibitem{Cao}
H. D. Cao, X. Sun, and Y. Xiaofeng, On the structure of gradient Yamabe solitons.
\textit{Math. Res. Lett.} \textbf{19} (2012),  767--774.


\bibitem{Chang&Chiu&Wu} S. C. Chang, H. L. Chiu, and C. T. Wu,
The Li-Yau-Hamilton inequality for Yamabe flow on a closed CR 3-manifold.
 \textit{Trans. Amer. Math. Soc.} \textbf{362} (2010), 1681--1698.

 \bibitem{Chang&Cheng} S. C. Chang and J. H. Cheng, The Harnack estimate for the Yamabe flow on CR manifolds
of dimension 3. \textit{Ann. Global Anal. Geom.} \textbf{21} (2002),
111--121.




\bibitem{Chow} B. Chow, The Yamabe flow on locally conformally flat manifolds
with positive Ricci curvature. \textit{Comm. Pure Appl. Math.}
\textbf{45} (1992), 1003--1014.

\bibitem{Chu} Y. Chu and X. Wang,
On the scalar curvature estimates for gradient Yamabe solitons. \textit{Kodai Math. J.} \textbf{36} (2013),  246--257.

\bibitem{Daskalopoulos} P. Daskalopoulos and N. Sesum,
 The classification of locally conformally flat Yamabe solitons. \textit{Adv. Math.} \textbf{240} (2013), 346--369.

\bibitem{DiCerbo} L. di Cerbo and M. Disconzi, Yamabe solitons, determinant of the Laplacian and the uniformization theorem for Riemann surfaces.
\textit{Lett. Math. Phys.} \textbf{83}  (2008), 13--18.

\bibitem{Gamara2} N. Gamara, The CR Yamabe conjecture---the case $n=1$. \textit{J.
Eur. Math. Soc.} \textbf{3} (2001), 105--137.


\bibitem{Gamara1} N. Gamara and R. Yacoub, CR Yamabe conjecture---the
conformally flat case. \textit{Pacific J. Math.} \textbf{201}
(2001), 121--175.

\bibitem{Ho1} P. T. Ho,
Result related to prescribing pseudo-Hermitian scalar
curvature. \textit{Int. J. Math.} \textbf{24} (2013), 29pp.

\bibitem{Ho2} P. T. Ho, The long time existence and convergence of the CR
Yamabe flow. \textit{Commun. Contemp. Math.} \textbf{14} (2012), 50 pp.

\bibitem{Hsu} S. Y. Hsu, A note on compact gradient Yamabe solitons. \textit{J. Math. Anal. Appl.}
 \textbf{388} (2012),  725--726.

\bibitem{Hsu1} S. Y. Hsu, Some properties of the Yamabe soliton and the related nonlinear elliptic equation.
\textit{Calc. Var. Partial Differential Equations} \textbf{49} (2014),  307--321.



\bibitem{Jerison&Lee1} D. Jerison and J.M. Lee, Extremals for the Sobolev inequality on
the Heisenberg group and the CR Yamabe problem. \textit{J. Amer.
Math. Soc.} \textbf{1} (1988), 1--13.

\bibitem{Jerison&Lee2} D. Jerison and J. M. Lee,
Intrinsic CR normal coordinates and the CR Yamabe problem.
\textit{J. Differential Geom.} \textbf{29} (1989), 303--343.

\bibitem{Jerison&Lee3} D. Jerison and J. M. Lee, The Yamabe problem on CR manifolds. \textit{J.
Differential Geom.} \textbf{25} (1987), 167--197.


\bibitem{Ma&Cheng}
L. Ma and L. Cheng, Properties of complete non-compact Yamabe solitons. \textit{Ann. Global Anal. Geom.} \textbf{40} (2011),  379--387.


\bibitem{Ma&Miquel} L. Ma and V. Miquel, Remarks on scalar curvature of Yamabe solitons. \textit{Ann. Global Anal. Geom.} \textbf{42} (2012),  195--205.

\bibitem{Schoen} R. Schoen, Conformal deformation of a Riemannian metric to constant scalar curvature.
\textit{J. Differential Geom.} \textbf{20} (1984), 479--495.

\bibitem{Schwetlick&Struwe} H. Schwetlick and M. Struwe, Convergence of the
Yamabe flow for ``large" energies.  \textit{J. Reine Angew. Math.}
\textbf{562} (2003), 59--100.



\bibitem{Trudinger} N. S. Trudinger, Remarks concerning the conformal deformation of Riemannian structures
on compact manifolds. \textit{Ann. Scuola Norm. Sup. Pisa (3)} \textbf{22} (1968), 265--274.

\bibitem{Ye} R. Ye, Global existence and convergence of Yamabe flow.
\textit{J. Diff. Geom.} \textbf{39} (1994), 35--50.




\end{thebibliography}

\end{document}